\documentclass{article}
\usepackage[utf8]{inputenc}
\usepackage{graphicx}
\def\MR#1{\href{http://www.ams.org/mathscinet-getitem?mr=#1}{MR#1}}

\usepackage{hyperref}
\usepackage{amsmath,amssymb,amsthm,mathrsfs}
\newtheorem{opr}{Definition}

\theoremstyle{plain}
\newtheorem{thm}{Theorem}
\newtheorem{lm}{Lemma}
\newtheorem{sld}{Corollary}
\newtheorem{utv}{Proposition}
\theoremstyle{plain}

\author{Emil Lerner\footnote{Affiliation: Moscow State University, Faculty of Computational Mathematics and Cybernetics, position: PhD student, email: neex.emil@gmail.com}
}
\title{The uniform distribution of sequences generated by iterated polynomials}
\date{}
\begin{document}
\maketitle
\begin{abstract}
Assume that $m,s\in\mathbb N$, $m>1$, while $f$ is a polynomial with integer coefficients, $\deg f>1$, $f^{(i)}$ is the $i$th iteration of the polynomial~$f$, $\kappa_n$ has a discrete uniform distribution on the set $\{0,1,\ldots,m^n - 1\}$. We are going to prove that with $n$ tending to infinity random vectors $\left(\frac{\kappa_n}{m^n},\frac{f(\kappa_n)  \bmod m^n}{m^n},\ldots,\frac{f^{(s - 1)}(\kappa_n)  \bmod m^n}{m^n}\right)$ weakly converge to a vector having a continuous uniform distribution in the $s$-dimensional unit cube. Analogous results were obtained earlier only for some classes of polynomials with $s\leqslant 3, \deg f = 2$.

The mentioned vectors represent sequential pseudorandom numbers produced by a polynomial congruential generator modulo $m^n$.
\end{abstract}

%\textbf{Mathematics Subject Classifications:} 05B35, 05C31, 11T06.


\section{Introduction}
\label{intro_section}
{We study the limit distribution of vectors of sequential pseudorandom numbers produced by a polynomial congruential generator modulo $m^n$ with the degree of the polynomial greater than $1$. The goal of this paper is to prove that with $n$ tending to infinity this limit distribution is uniform in the $s$-dimensional cube (with any~$s$). As appeared, this property takes place with any $m$ (in applications we put $m=2$).}

Let $m$ be a fixed positive integer, $m>1$, $[m]=\{0,\ldots,m-1 \}$. For an integer number~$x$ we denote by $x\bmod m$ the least nonnegative residue of $x$ modulo $m$.
\begin{opr}
A function $f:\mathbb Z\to\mathbb Z$ is said to be \textit{compatible} if for any $n\in\mathbb N$ and $x,y\in \mathbb Z$ the equality $x\bmod m^n=y\bmod m^n$ implies that $f(x)\bmod m^n=f(y)\bmod m^n$.
\end{opr}

\noindent
Below in this section we consider only compatible functions~$f$. Evident examples of compatible functions are polynomials with integer coefficients.

Denote iterations of a function $f$ by $f^{(i)}$, (i.e., $f^{(0)}(x) \equiv x$, $f^{(1)}(x) \equiv f(x)$, $f^{(2)}(x) \equiv f(f(x))$, etc). Note that for a compatible function~$f$ any iteration $f^{(i)}$ also represents a compatible function.

For any compatible function~$f$ we put $$\mathbf P_n^s(f)=\left\{
\left(\frac{x}{m^n},\frac{f(x) \bmod m^n}{m^n},\ldots,\frac{f^{(s-1)}(x) \bmod m^n}{m^n}\right), x\in [m^n] \right\}. $$ The set $$\mathbf P^s(f) = \bigcup\limits_{n=1}^\infty \mathbf P_n^s(f)$$ is called the $s$-dimensional projection of the function $f$.

In \cite{anashin_01law}, \cite{Anashin} one proves the following important theorem.

\begin{thm}[rule $0-1$]\label{zeroonelaw}
Let $m$ be a prime number. For any compatible function $f$ the Lebesgue measure of the closure of its two-dimensional projection equals either~0 or~1.
\end{thm}

It is well known \cite[section~3.3.4]{knuth} that if $f$ is a polynomial of degree~1, then the Lebesgue measure of the closure of its projection equals~0, and with $n$ tending to infinity all points appear to be located in several hyperplanes inside the unit hypercube. Moreover, in \cite{Anashin} it is also proved that 
{ with $\deg f>1$} the measure of the closure of the two-dimensional projection of each polynomial with integer coefficients equals 1. We are interested in a more difficult question, namely, we study the distribution of vectors
\begin{equation}
\label{vecfirst}
\left(\frac{x}{m^n},\frac{f(x) \bmod m^n}{m^n},\ldots,\frac{f^{(s-1)}(x) \bmod m^n}{m^n}\right)
\end{equation}
with $x$ randomly chosen in the set~$[m^n]$.

Recall that for a random vector $\eta=(\eta(1),\ldots,\eta(s))$ the value of the cumulated
{distribution} function (CDF) $F(x_1,\ldots,x_n)$ equals the probability of the random event
$${\mathbb P}(\,\eta(1)\leqslant x_1,\ldots,\eta(s)\leqslant x_s).$$
A continuous uniform distribution in the~$s$-dimensional unit cube obeys the CDF $F(x_1,\ldots,x_n)=\prod_{i=1}^n U(x_i)$, where
$$
U(x_i)=\left\{\begin{matrix}
0, &\mbox{ if } x_i<0, \\
x_i, &\mbox{ if } x_i \in [0,1],\\
1, &\mbox{ if } x_i>1.
\end{matrix}
\right.
$$

A discrete uniform distribution on a finite set~$\mathbf X\subset \mathbb R^s$ is defined by equal probabilities of all values in this set, i.e., ${\mathbb P}(\mathbf x)=1/|\mathbf X|$, $\mathbf x\in \mathbf X$. The CDF of a discrete uniform distribution is stepwise~\cite{doob}. In what follows, we essentially use the finiteness of the set of values of a discrete uniform random variable (r.v.), as distinct from the sets of values of a continuous uniform r.v. However, in this paper we meet a continuous uniform distribution (in the $s$-dimensional unit cube) only once, namely, when proceeding to the limit in a discrete case. Recall that a weak convergence of $s$-dimensional random vectors $\eta_n$ to some random vector $\eta$ means that CDFs of $\eta_n$ pointwisely converge to the CDF of~$\eta$ for every point $\mathbf x\in \mathbb R^s$ at which CDF of~$\eta$ is continuous.
Note that the weak convergence property of CDFs of considered random vectors is independent of the probabilistic space, in which these vectors are given.

\begin{thm}[The main theorem]
\label{mainthm0}
For any polynomial $f$ with integer coefficients whose degree is greater than 1 and for any positive integer $s$ the sequence of random vectors $$\left(\frac{\kappa_n}{m^n},\frac{f (\kappa_n)\bmod m^n}{m^n},\ldots,
\frac{f^{(s-1)}(\kappa_n)\bmod m^n}{m^n}\right),$$ where $\kappa_n$ are discrete uniform
distributions on the finite set $[m^n]$, {weakly converges (as $n$ tends to infinity) to the continuous uniform distribution on~$[0, 1)^s$.}
\end{thm}

More strong results are obtained (in other terms) in papers~\cite{nemcy}, \cite{nemcy2} as particular cases of Theorem~\ref{mainthm0}; they deal with second-degree polynomials with $s \in \{2, 3\}$. In~\cite{blastra} these results are improved for the case of two iterations and $m=2$. This work is first to study the general case.

The paper has the following structure. In Section~\ref{new_section} we introduce a probabilistic space which is 
{convenient for the proof of Theorem~\ref{mainthm0}} and propose a new statement of {this} theorem. In Section~\ref{basicdef_section} we generalize the introduced notions for the case of an arbitrary collection of compatible functions (instead of iterations of one and the same function). The key moment in our proof is an analog of Theorem~\ref{mainthm0} for the collection of monomials $x, x^2, \ldots, x^{s}$; this result is proved in Section~\ref{monoms_section}. In three previous sections, we study the necessary and sufficient conditions for the weak convergence of an arbitrary collection of compatible functions to a continuous uniform distribution. In Section~\ref{polynoms_section}, using previous results, we easily prove the weak convergence for an arbitrary collection of polynomials of various degrees and, as a corollary, get the assertion of Theorem~\ref{mainthm0}.

\section{The statement of the main theorem in terms of $m$-adic numbers}
\label{new_section}
Recall that any integer $m$-adic number $x$ can be written as an infinite sequence of digits that belong to the set $[m]=\{0,\ldots,m-1 \}$ (digits of a ``left-infinite natural number''), operations of addition and column multiplication in $\mathbb Z_m$ can be defined analogously to usual operations of addition and column multiplication in the $m$-ary notation. The ring~$\mathbb Z$ considered in the previous section is a part of~$\mathbb Z_m$ (see, for example,~\cite{mahler} for more detail).

Denote the set of all $m$-adic numbers by $\mathbb Q_m$; such numbers are representable as formal power series
\begin{equation}
\label{series}
x=\sum_{i=k}^\infty a_i m^i,
\end{equation}
where $a_i\in [m]$, $k\in\mathbb Z$ (see~\cite{Gelfand}). If, in addition, $x\in{\mathbb Z}_m$, then $k\geqslant 0$.

If $m$ is prime, then $\mathbb Q_m$ is a field (in this case, one writes the symbol~$p$ in place of $m$). In a general case, $\mathbb Q_m$, as well as $\mathbb Z_m$, represents a ring~\cite{mahler}.

The ring $\mathbb Z_m$ is a metric space whose metric obeys the formula $\rho(x,y)=m^{-n}$, where~$n$ is the first position, at which integer $m$-adic numbers~$x$ and~$y$ start to differ. Recall that~(\cite{mahler},~\cite{Gelfand}) for any $x\in\mathbb Z_m$ the number $x\bmod m^n$ is defined as the only number $y\in [m^n]$ such that $\rho(x,y)\leqslant m^{-n}$. In other words, if~$x$ is written in form~(\ref{series}) with $k=0$, then $y=\sum_{i=0}^{n-1} a_i m^i$. This definition corresponds to the definition of the function $\bmod\, m^n$ as a map from $\mathbb Z$ to $[m^n]$ (see \cite{mahler},~\cite{Gelfand}); recall that the latter definition is given in the previous section.

\begin{opr}
A function $f:$  ${\mathbb Z}_m\to{\mathbb Z}_m$ is said to be \textit{1-Lipschitz} if for any $\varepsilon>0$ the inequality $\rho(x,y)<\varepsilon$ gives $\rho(f(x),f(y))<\varepsilon$.
\end{opr}

Therefore, the 1-Lipschitz property of a function $f$ means that given $n$ junior digits of some number $x$, one can uniquely define $n$ junior digits of the number $y=f(x)$. Any function~$f$ with such a property generates a collection of functions $(f \bmod m^n)$ acting from $[m^n]$ to $[m^n]$ which obey the formula $(f \bmod m^n): x \bmod m^n\to y\bmod m^n$. Note that for $x \in [m^n]$ it holds that $(f \bmod m^n)(x) = f(x) \bmod m^n$.

Evidently, notions of the compatibility and 1-Lipschitz property are interconnected. The next proposition describes this interconnection. 

\begin{utv}[\cite{mahler}]
Any compatible function acting from~$\mathbb Z$ to $\mathbb Z$ is uniquely extendable to an 1-Lipschitz function acting from ${\mathbb Z}_m$ to ${\mathbb Z}_m$.
\end{utv}
\par\noindent\textbf{Proof:}
Any 1-Lipschitz function is continuous, therefore it suffices to define it on an everywhere dense set like the set $\mathbb Z$ in $\mathbb Z_m$. By definition, any compatible function satisfies 1-Lipschitz conditions on the set $\mathbb Z$.
$\square$

In particular, polynomials with integer coefficients are extendable to 1-Lipschitz functions on the set of integer $m$-adic numbers. Denote the extension of a compatible function $f$ acting from $\mathbb Z$ to $\mathbb Z$ up to an 1-Lipschitz function acting from ${\mathbb Z}_m$ to ${\mathbb Z}_m$ by the same symbol~$f$. In what follows we consider only {1-Lipschitz} functions $f$ acting from $\mathbb Z_m$ to $\mathbb Z_m$
{(}unless otherwise is specified{). Moreover, we assume that $f$ extends a compatible function acting from $\mathbb Z$ to $\mathbb Z$.} 

Let us give a new statement of the main theorem. Note that with $x\in [m^n]$ for any 1-Lipschitz function~$f$ vector~(\ref{vecfirst}) by definition coincides with the vector
\begin{equation}
\label{vecsecond}
\left(\frac{x}{m^n},\frac{(f \bmod m^n)(x)}{m^n},\ldots,\frac{(f^{(s-1)}\bmod m^n)(x)}{m^n}\right).
\end{equation}

Assume that $\Omega$ is the space of elementary events~$\mathbb Z_m$ whose probabilistic measure equals the standard Haar measure $\mu$ normalized by the condition $\mu(\mathbb Z_m)=1$ (see~\cite{Gelfand}). Recall that measurable sets in this probabilistic space are countable unions of balls $a+m^k \mathbb Z_m$, where $k\in\mathbb N$, $a\in \mathbb Z_m$, while $\mu(a+m^k \mathbb Z_m)=m^{-k}$.

Let us define random values $\xi_n$ by the formula~$\xi_n(\omega)=\omega \bmod m^n$; here the symbol $\omega$ denotes an elementary event, i.e., $\omega\in \mathbb Z_m$. The probability $\mathbb P(\xi_n = a)$ for $a \in [m^n]$ is the probabilistic measure of the set defined by the condition $\omega \bmod m^n =a$, i.e., it is the standard Haar measure of the ball $a + \mathbb Z_m m^n$. This measure is independent of $a$ (and equals $m^{-n}$), therefore random values~$\xi_n$ have a discrete uniform distribution on the finite set $[m^n]$. Moreover, they are connected by correlations 
$\xi_{n} = \xi_{n+1} \bmod m^n$, $n\in\,\mathbb N$.

\begin{thm}
\label{mainthm}
For any polynomial $f$ with integer coefficients of a degree greater than~1 and for any positive integer $s$ the sequence of random vectors
\begin{equation}
\label{vec}
\left(\frac{\xi_n}{m^n},\frac{(f \bmod m^n)(\xi_n)}{m^n},\ldots,\frac{(f^{(s-1)}\bmod m^n)(\xi_n)}{m^n}\right)
\end{equation}
weakly converges (as $n$ tends to infinity) to a random vector with a continuous uniform distribution on $[0, 1)^s$.
\end{thm}

The next proposition follows from that proved above.
\begin{utv}
Theorems~\ref{mainthm} and~\ref{mainthm0} are equivalent. \label{ochevidno2}
\end{utv}

\par\noindent\textbf{Proof:}
As was mentioned earlier, for the polynomial~$f$ which is studied in theorems~\ref{mainthm} and~\ref{mainthm0} formulas~(\ref{vecfirst}) and~(\ref{vecsecond}) are equivalent with any $x\in [m^n]$. Random values $\xi_n$ and $\kappa_n$ are distributed identically, namely, their distribution is discrete uniform on the finite set $[m^n]$. The notion of a weak convergence depends only on distributions of random values (it is independent of the probabilistic space, where these values are defined).
$\square$

We are going to prove the main theorem just in the proposed statement. The definition of elementary events on~$\mathbb Z_m$ allows us to effectively estimate the measure of events, for which the $i$th coordinate of vector~(\ref{vec}) does not exceed $z_i$, where $(z_1,\ldots,z_s)$ is an arbitrary vector in $[0,1)^s$. The statement of Theorem~\ref{mainthm} means that this measure (with $n$ tending to infinity) converges to $z_1\times \ldots\times z_s$, i.e., to $V(J)$, where the symbol $J$ stands for the parallelepiped $\{(x_1,\ldots,x_s): x_i\in [0,z_i], i=1,\ldots,s\}$, and $V(J)$ does for its $s$-dimensional volume.

\section{The collection of functions that are uniform with each suffix}
\label{basicdef_section}
In this paper, instead of iterations of a certain function $f$, we often consider an arbitrary collection of functions $f_1,f_2,\ldots,f_s$. We study random vectors in~$[0,1)^s$ in the form
\[
\varphi_n (f_1,\ldots,f_s)=\left(\frac{(f_1 \bmod m^n)(\xi_n)}{m^n},\ldots,\frac{(f_s\bmod m^n)(\xi_n)}{m^n}\right),
\]
where random values $\xi_n$ are defined at the end of the previous section. Note that in view of the 1-Lipschitz property of functions~$f_i$ the vector~$\varphi_n (f_1,\ldots,f_s)$ is representable as a function of an elementary event~$\omega\in \mathbb Z_m$, namely,
\begin{equation}
\label{otherreptresentation}
\varphi_n (f_1,\ldots,f_s)=\left(\frac{f_1(\omega) \bmod m^n}{m^n},\ldots,\frac{f_s(\omega)\bmod m^n}{m^n}\right).
\end{equation}

\begin{opr} A collection of functions $(f_1,\ldots,f_s)$ is called a {\it generating uniform distribution}, if with $n$ tending to infinity $\varphi_n(f_1,\ldots,f_s)$ weakly converges to a random vector with a continuous uniform distribution in $[0,1)^s$. \end{opr}

Let the symbol $\mathfrak J_s$ denotes the set of all parallelepipeds such that they belong to the cube $[0, 1)^s$ and their sides are parallel to coordinate axes; denote the $s$-dimensional volume of a parallelepiped $J \in \mathfrak J_s$ by the symbol $V(J)$.

\begin{utv}
A collection of functions $(f_1,\ldots,f_s)$ {\it generates a uniform distribution} if and only if for any parallelepiped $J \in \mathfrak J_s$,
\begin{equation}
\label{form}
\lim_{n\to\infty} {\mathbb P}(\varphi_n(f_1,\ldots,f_s)\in J) = V(J).
\end{equation}\label{utv_prosto}
\end{utv}
\par\noindent\textbf{Proof:} A weak convergence of $\varphi_n(f_1,\ldots,f_s)$ to a vector having a continuous uniform distribution in $[0,1)^s$ represents a particular case of formula~(\ref{form}) for parallelepipeds \begin{equation}\label{vida}\{(x_1,\ldots,x_s): x_i\in [0,z_i], i=1,\ldots,s\},
\end{equation} 
where $(z_1, \ldots, z_s) \in [0, 1]^s$. Therefore, formula~(\ref{form}) implies a weak convergence. Vice versa, any parallelepiped in the set $\mathfrak J_s$ belongs to the algebra of sets generated by parallelepipeds~(\ref{vida}).$\square$

Let $k$ be some fixed number, $a_1,a_2,\ldots,a_s \in [m^k]$, $\mathbf a = (a_1, a_2, \ldots, a_s)$. Denote by $J_k(\mathbf a)$ the cube in $[0,1)^s$ defined by inequalities
\begin{eqnarray*}
\frac{a_1}{m^k} &\leqslant z_1 <& \frac{a_1+1}{m^k},\\
\frac{a_2}{m^k} &\leqslant z_2 <& \frac{a_2+1}{m^k},\\
&\ldots&\\
\frac{a_s}{m^k} &\leqslant z_s <& \frac{a_s+1}{m^k};
\end{eqnarray*}
here $(z_1,z_2,\ldots,z_s)$ are coordinates of a point in $[0,1)^s$.

Since any parallelepiped $J \in \mathfrak J_s$ can be arbitrarily accurate approximated by the union of several parallelepipeds in form $J_k(\mathbf a)$, Proposition~\ref{utv_prosto} implies the following one.

\begin{utv} \label{opr_orig} Fix a collection $f_1, f_2, \ldots, f_s$. If for all $J = J_k(\mathbf a)$ for any $k \in \mathbb N, \mathbf a \in [m^k]^s$,
\begin{equation}
\label{maindef}
\lim_{n\to\infty} {\mathbb P}(\varphi_n(f_1,\ldots,f_s)\in J) = V(J) 
\end{equation}
then the collection $f_1, f_2, \ldots, f_s$ generates a uniform distribution.
\end{utv}
Let $x\in\mathbb Z_m$, $i,j\in\mathbb N$, $i>j$. Define $substr(x,i,j)$ as
$$
substr(x,i,j)=\frac{x\bmod m^i-x\bmod m^j}{m^j}.
$$
In other words, the number $substr(x,i,j)$ belongs to $[m^{i-j}]$ and is formed by digits that correspond to degrees of~$m$ varying from $j$ to $i-1$. Note that $substr(x,i,j)$ is also defined for $x\in \mathbb N_0\subset\mathbb Z_m$, in particular, for $x\in[m^i]$.

Let us generalize the introduced denotation for vectors $\mathbf x=(x_1,\ldots,x_s)$,  ${x_1,\ldots,x_s\in\mathbb Z_m}.$ Put
$$
substr(\mathbf x,i,j)=(substr(x_1,i,j),\ldots,substr(x_s,i,j)).
$$
In particular, we can write the definition of $\xi_n$ as $\xi_n= substr(\omega,n,0)$.

Let $n,k\in \mathbb{N},n \geqslant k, \mathbf x \in [m^n]^s, \mathbf a \in [m^k]^s$. Note that if a point $$\left(\frac{x_1}{m^n},\frac{x_2}{m^n},\ldots,\frac{x_s}{m^n}\right)$$ belongs to $J_k(\mathbf a)$, then 
$
substr(\mathbf x,n,n-k)=\mathbf a.
$

\label{basic_def}

In what follows, for convenience, instead of unconditional probabilities 
$${\mathbb P}(\varphi_n(f_1,\ldots,f_s)\in J),$$ we deal with conditional ones ${\mathbb P}(\varphi_n(f_1,\ldots,f_s)\in J\,|\, \xi_d=\beta)$, where $\beta\in [m^d]$ for some positive integer~$d$. Evidently, the condition $x \bmod m^d=\beta$ is equivalent to the coincidence of $substr(x,d,0)$ (the suffix in the $m$-ary notation) with $\beta$.

\begin{opr} We treat the collection of 1-Lipschitz functions $f_1,\ldots,f_s$ as a {\it generating uniform distribution with each suffix}, or just a {\it uniform distribution with each suffix}, if \begin{equation}
\label{mainmaindef}
\lim_{n\to\infty}{\mathbb P}(\varphi_n(f_1,\ldots,f_s)\in J\,|\, \xi_d=\beta) = V(J)
\end{equation}
for any $J = J_k(\mathbf a)$, where $k \in \mathbb N, \mathbf a \in [m^k]^s$, and any~$d \in \mathbb N, \beta\in[m^d]$.
\end{opr}

\begin{utv} If a collection of 1-Lipschitz functions $f_1,\ldots,f_s$ is uniform with each suffix, then it generates a uniform distribution. \label{utv_ochevidno}
\end{utv}
\par\noindent\textbf{Proof:} Let us consider a complete group of events $\xi_d = \beta$, where $\beta \in [m^d]$ with some fixed positive integer $d$. By the total probability formula,
$${\mathbb P}(\varphi_n(f_1,\ldots,f_s)\in J) = \sum_{\beta \in [m^d]} {\mathbb P}(\varphi_n(f_1,\ldots,f_s)\in J\,|\, \xi_d=\beta)/m^d.$$ Therefore from~(\ref{mainmaindef}) we get~(\ref{maindef}).
$\square$

Let us state the key theorem of this paper.

\begin{thm}[The key theorem]
\label{monoms}
Let $s \in \mathbb{N}$. The collection of functions
$$
f_1(x)=x,\ f_2(x)=x^2,\ldots,f_s(x)=x^s
$$
is uniform with each suffix.
\end{thm}

\section{Corollaries of the uniformity with each suffix}
\label{coin_section}
Let us prove Theorem~\ref{monoms} by induction with respect to~$s$. For performing the inductive transition, we need some corollaries of the uniformity.

Let us first state a simple probabilistic proposition.

\begin{lm}
\label{lem:Levy}
Let an infinite collection of events $A_1,A_{2},\ldots$ be such that there exists $\varepsilon > 0$ such that for any $t\in\mathbb N$,
\begin{equation}
\label{eqlemLevy}
\overline{\lim}_{n\to\infty} {\mathbb P}
\left( A_n\,|\,\overline{A_1},\ldots,\overline{A_t}\,\right)>\varepsilon,
\end{equation}
(here the symbol $\overline{\lim}$ denotes the upper limit, and $\overline{A}$ does the event opposite to~$A$). Then the probability of the event $\bigcup_{i=1}^n A_i$ (note that it takes place if and only if so does at least one of events from the collection $A_1,\ldots,A_n$) tends to 1 as $n\to\infty$.
\end{lm}

\par\noindent\textbf{Proof:}
The event that is opposite to~$\bigcup_{i=1}^n A_i$ is $\bigcap_{i=1}^n \overline{A_i}$. Its probability, evidently, decreases with the increase of~$n$, so for proving the lemma it suffices to make sure that this probability cannot be bounded from below. In particular, it suffices to find (for arbitrary~$t$) a number~$n$ such that
$${\mathbb P}\left(\bigcap\nolimits_{i=1}^n \overline{A_i}\,\right)< {\mathbb P}\left(\bigcap\nolimits_{i=1}^t \overline{A_i}\,\right)(1-\varepsilon).$$

By the condition of the lemma,
$$
\underline{\lim}_{\ n\to\infty} {\mathbb P}
(\overline{A_n}\,|\,\overline{A_1},\ldots,\overline{A_t}\,)< 1-\varepsilon.
$$
By the definition of the notion of conditional probability, this means that there exists~$n$ such that
$${\mathbb P}\left(\left(\bigcap\nolimits_{i=1}^t \overline{A_i}\right) \cap \overline{A_n}\right)< {\mathbb P}\left(\bigcap\nolimits_{i=1}^t \overline{A_i}\right)(1-\varepsilon).$$
The left-hand side of the latter equality, evidently, gives an upper bound of the desired probability
${\mathbb P}\left(\bigcap_{i=1}^n \overline{A_i}\,\right)$.
$\square$

Note that the proved lemma is a certain simplified version of the Levy theorem (see \cite[section 7.4]{doob}, theorem~4.1, corollary~2).

\begin{sld}
\label{cond:Levy}
Let $B$ be some fixed event, ${\mathbb P}(B)>0$, and $j\in\mathbb N$. Let an infinite collection of events $C_j,C_{j+1},\ldots$ be such that there exists $\varepsilon' > 0$, with which for any $t\in\mathbb N, t\geqslant j$,
\begin{equation}
\label{cond_sld}
\lim_{n\to\infty} {\mathbb P}
\left( C_n\,|\,\overline{C_j},\ldots,\overline{C_t},\,B\right)=\varepsilon',
\end{equation}
Then
$
\lim_{n\to\infty} {\mathbb P}
\left(\bigcup\nolimits_{i=j}^n C_i\,|\,B\,\right)= 1.
$
\end{sld}

\par\noindent\textbf{Proof:}
Let us define (on the same sigma-algebra of events, were the measure $\mathbb P$ is given) a new probabilistic measure ${\mathbb P}'$ by the conditional probability formula. To this end, for any event~$D$ let us put ${\mathbb P}'(D)={\mathbb P}(D\cap B)/{\mathbb P}(B)$. The sigma-additivity of the new probabilistic measure follows from the sigma-additivity of the initial measure~${\mathbb P}$, because we can take the constant $1/{\mathbb P}(B)$ out of the sum sign; the property of the probabilistic normalization of the measure ${\mathbb P}'$ for the whole space of elementary events $\Omega$, i.e., ${\mathbb P}'(\Omega)=1$, is evident by definition.

Let us apply Lemma~\ref{lem:Levy} to this new probabilistic space, considering events $C_j, C_{j+1},\ldots$ in place of those~$A_1,A_2,\ldots$. By definition,
$$
{\mathbb P}\left( C_n\,|\,\overline{C_j},\ldots,\overline{C_t},\,B\right)=
{\mathbb P}'\left( C_n\,|\,\overline{C_j},\ldots,\overline{C_t}\right),
$$
and, consequently, inequalities~(\ref{eqlemLevy}) are fulfilled for $\varepsilon \in (0, \varepsilon')$. We get
$$
{\mathbb P}
\left(\bigcup\nolimits_{i=j}^n C_i\,|\,B\,\right)={\mathbb P}'
\left(\bigcup\nolimits_{i=j}^n C_i\,\right),
$$
i.e., the desired assertion follows from the proposition of the lemma in the considered case.
$\square$

Let apply this corollary to the uniformity property with each suffix.

\begin{lm}
\label{flipflop_forone}
Let a collection of functions $f_1,f_2, \ldots,f_s$ generate a uniform distribution with each suffix. Choose arbitrary $\beta\in[m^d]$, $\mathbf a \in [m^k]^s$, where $d,k \in \mathbb N$. Then for any $j\in \mathbb N$ the probability
\begin{equation}
\label{limP}
{\mathbb P}(\,\bigcup\nolimits_{i=j}^n (\varphi_i(f_1,\ldots,f_s)\in J_k(\mathbf a))\,|\,\xi_d=\beta)
\end{equation}
tends to 1 as $n$ tends to infinity.
\end{lm}

\par\noindent\textbf{Proof:}
Let us make use of Corollary~\ref{cond:Levy}. Consider the event $\xi_d=\beta$ for $B$ and do events $\varphi_i(f_1,\ldots,f_s)\in J_k(\mathbf a)$ for $C_i$. Note that if for some positive integer~$t$, $t\geqslant j$, it holds that ${\mathbb P} \left( \left(\bigcap\nolimits_{i=j}^t \overline{C}_i\right)\cap B\right)=0$, then the unit limit value of the nondecreasing probability~(\ref{limP}) is attained as early as at~${n=t}$. Let us make sure that otherwise conditions~(\ref{cond_sld}) are fulfilled.

Really, in view of the 1-Lipschitz property of functions $f_1,\ldots,f_s$ and the definition of $\varphi_i$, the considered nonempty event $\left(\bigcap\nolimits_{i=j}^t \overline{C}_i\right)\cap B$ is representable in terms of random values $\xi_t$; more precisely, it is representable as the union of events $\bigcup_{\beta'\in\mathcal B}\{\omega:\xi_t=\beta'\}$ for some $\mathcal B\subseteq [m^t]$ (here we assume that $t\geqslant d$, one can treat the case of the opposite inequality just in the same way).

By condition (the generation of the uniform distribution with each suffix) with any $\beta'\in [m^t]$,
\begin{equation}
\label{lim}
\lim_{n\to\infty} {\mathbb P}
\left( \varphi_n(f_1,\ldots,f_s)\in J_k(\mathbf a)\,|\,\xi_t=\beta'\right)=m^{-ks}.
\end{equation}
Evidently, for noncoinciding $\beta'$ events $\xi_t=\beta'$ are incompatible. Summing up equalities~(\ref{lim}) for all $\beta'\in \mathcal B$ and dividing by $|\mathcal B|$, we conclude that for any $\mathcal B\subseteq [m^t]$, 
$$
\lim_{n\to\infty} {\mathbb P}
\left( \varphi_n(f_1,\ldots,f_s)\in J_k(\mathbf a)\,|\,\xi_t\in \mathcal B\right)=m^{-ks}.
$$
Therefore, conditions~(\ref{cond_sld}) with $\varepsilon'=m^{-ks}$ are fulfilled.
$\square$

For convenience of further considerations let us state Lemma~\ref{flipflop_forone} in a different form.

\begin{sld}
\label{sldcond}
Let a collection of functions $f_1,f_2, \ldots,f_s$ generate a uniform distribution with each suffix. Fix arbitrarily $j,d,k \in \mathbb N$, $\beta\in[m^d]$, $\mathbf a \in [m^k]^s$, and a small value $\varepsilon$. Then there exists $N$ such that 
\begin{eqnarray}
\nonumber
&\mbox{among $x\in [m^N]$, $x\bmod m^d=\beta$, the relative amount of~$x$ such that}\\
&\label{cond} \mbox{$substr((f_1(x),\ldots,f_s(x)),\ell,\ell-k)=\mathbf a$ with some $\ell$, $j\leqslant \ell\leqslant N$,}\\
\nonumber
&\mbox{exceeds $1-\varepsilon$.}
\end{eqnarray}
\end{sld}

Note that we need the notion of the uniformity with each suffix just for obtaining results described in this section. Results described in the following sections are also valid for a collection of functions that generate a uniform distribution (not necessarily with each suffix). However, we consider them mainly for the case of the uniformity with each suffix; this allows us to prove Theorem~\ref{mainthm}.

\section{Auxiliary results for the inaccuracy of hitting the cube $J_k(\mathbf a)$}
\label{cubic_section}
Theorem~\ref{monoms} means that the collection of functions $(x,x^2,\ldots,x^s)$ satisfies correlation~(\ref{mainmaindef}). In Lemma~\ref{flipflop_forone} we consider a corollary of this correlation; we use it in the inductive passage with respect to~$s$, assuming (as the induction hypothesis) that it is valid for the collection $(x,x^2,\ldots,x^{s-1})$. However here we consider sufficient conditions, whose fulfillment for the whole collection of~$s$ monomials is to be proved.

The problem is that the ``aimed hitting'' of senior positions in the cube $J_k(\mathbf a)$ is complicated because of carryovers caused by the ``accumulation'' that takes place in junior positions. In this section we first prove a simple auxiliary assertion (Lemma~\ref{stolbik}), which restricts this inaccuracy, and then state sufficient conditions for~(\ref{mainmaindef}) under constraints imposed on the ``hitting inaccuracy''.

The assertion of Lemma~\ref{stolbik} restricts the variation of senior positions in the summation or subtraction of two numbers. Here we understand the variation as the minimum of two differences modulo $m^k$, i.e., we say that 0 and $m^k-1$ differ by 1.

\begin{lm}
\label{stolbik}
Let $n, k \in \mathbb N, n \geqslant k, x_1, x_2 \in \mathbb Z$. Put $y_1 = substr(x_1, n, n-k), y_2 = substr(x_2, n, n-k)$. Then we conclude that \begin{eqnarray}
&&(substr(x_1+x_2, n, n-k) - (y_1 + y_2)) \bmod m^k \in \{0, 1\};\label{st1}\\
&&(substr(x_1-x_2, n, n-k) - (y_1 - y_2)) \bmod m^k \in \{0, m^k-1\}.\label{st2}
\end{eqnarray}
\end{lm}
\par\noindent\textbf{Proof:}
{Since $substr(x_1\pm x_2, n, n-k)=substr(x_1\bmod m^n \pm x_2\bmod m^n, n, n-k)$, it suffices to consider the case of $x_1,x_2\in [m^n]$.}
Both {desired} propositions follow from procedures of column summation and subtraction of numbers, when the carryover to senior $k$ positions in summation (as well as the borrowing in subtraction) does not exceed 1. Therefore, senior $k$ positions in the sum (difference) differ from the sum (difference) of senior $k$ positions modulo $m^k$ at most by 1.
$\square$

Let us now return to sufficient conditions for generating a uniform distribution. Let $K\in\mathbb N$, $\mathbf b \in [m^K]^s$. Denote by $O_K(\mathbf b)$ a neighborhood of a vector $\mathbf b$, more precisely, the totality of all collections $\mathbf c\in [m^K]^s$ such that $(b_i-c_i)\bmod m^K \in \{0, 1, m^K-1\}$ for all $i=1,\ldots,s$. For example, in these terms, correlations~(\ref{st1}),~(\ref{st2}) imply that
$y_1\pm y_2\in O_k(substr(x_1\pm x_2, n, n-k))$.

\begin{lm}
\label{OK}
A sufficient condition for the uniformity of a collection of functions $f_1,\ldots,f_s$ with each suffix is that for any positive integers $K$ and $d$, $\beta\in [m^d]$, and for any $M\subseteq [m^K]^s$, 
\begin{equation}
\label{keyeneq}
\underline\lim_{\,n\to\infty}{\mathbb P}(\varphi_n(f_1,\ldots,f_s)
\in \bigcup_{\mathbf b\in M} \bigcup_{\mathbf c\in O_K(\mathbf b)}J_K(\mathbf c)\,|\, \xi_d=\beta)\geqslant
\sum_{\mathbf b\in M} V(J_K(\mathbf b)).
\end{equation}
\end{lm}
\par\noindent\textbf{Proof:}
Let us prove correlation~(\ref{mainmaindef}) for fixed $k\in\mathbb N, \mathbf a\in [m^k]^s, d \in \mathbb N, \beta \in [m^d]$. For brevity, let us introduce the denotation
$$
p_n(\mathbf a)={\mathbb P}(\varphi_n(f_1,\ldots,f_s)\in J_k(\mathbf a)\,|\, \xi_d=\beta).
$$

\begin{figure}[bh]
\noindent\centering{
\includegraphics[width=50mm]{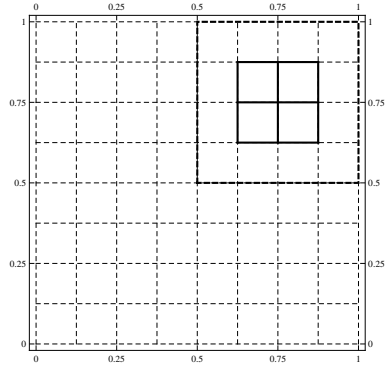}
}
\caption{Visualization of $M(K)$ for the case of $m=2$, $s=2$, $k=1$, $\mathbf a=(1,1)$, and $K=3$. Dashed bold lines and solid bold lines are used, respectively, for highlighting the square $J_k(\mathbf a)$ and squares $J_K(\mathbf b)$ for all $\mathbf\ b\in M(K)$.
}
\label{pic1}
\end{figure}

Let $K>k$. Denote by $M(K)$ the totality of vectors $\mathbf b$, $\mathbf b\in [m^K]^s$, such that $J_K(\mathbf c)\subseteq J_k(\mathbf a)$ for any $\mathbf c\in O_K(\mathbf b)$ (see Fig.~\ref{pic1}). Evidently, $M(K)\subseteq M'(K)$, where $M'(K)$ is the set of all vectors $\mathbf b$ from $[m^K]^s$ such that $substr(\mathbf b,K,K-k)=\mathbf a$, or, equivalently, $J_K(\mathbf b)\subseteq J_k(\mathbf a)$. Certainly, $|M'(K)|=m^{(K-k)s}$.

Let us now calculate $|M(K)|$. The set $M(K)$ consists of vectors $\mathbf b=(b_1,\ldots,b_s)$ (and only of them) such that $substr(\mathbf b, K, K-k) = \mathbf a$ and $$b_i \bmod m^{K-k} \in \{1, \ldots, m^{K-k}-2\}$$ for all~$i, 1 \leqslant i \leqslant s$. Since the cardinal number of this set equals $m^{K-k}-2$, we conclude that $|M(K)|=(m^{K-k}-2)^s$.

Therefore, $|M'(K)\setminus M(K)| = m^{(K-k)s} - (m^{K-k}-2)^s$. With $K$ tending to infinity this value is $o(m^{(K-k)s})=o(m^{Ks})$.

Since $M(K) \subseteq M'(K)$, condition~(\ref{keyeneq}) implies that for any~$K$,
\begin{equation}
\label{sledkeyneq}
\underline\lim_{\,n\to\infty}p_n(\mathbf a)\geqslant
\sum_{\mathbf b\in M(K)} V(J_K(\mathbf b)).
\end{equation}
Thus,
$$
V(J_k(\mathbf a))-\sum_{\mathbf b\in M(K)} V(J_K(\mathbf b))=
\sum_{\mathbf b\in M'(K)\setminus M(K)} V(J_K(\mathbf b)).
$$
With $K$ tending to infinity the right-hand side is $o(m^{Ks})m^{-Ks} = o(1)$. Therefore, from~(\ref{sledkeyneq}) we get the inequality
$$
\underline\lim_{\,n\to\infty}p_n(\mathbf a)\geqslant
V(J_k(\mathbf a)).
$$

Furthermore, since $[0, 1)^s$ is the union of cubes $J_k(\mathbf a)$, $\mathbf a\in [m^k]^s$, we get
$$
 p_n(\mathbf a)=1-
\sum_{\mathbf b\in [m^k]^s \setminus\{\mathbf a\}}
 p_n(\mathbf  b),\qquad V(J_k(\mathbf a))=1-\sum_{\mathbf b\in [m^k]^s\setminus\{\mathbf a\}} V(J_k(\mathbf b)).
$$
Consequently,
$$
\overline\lim_{\,n\to\infty} p_n(\mathbf a)\leqslant 1-\!\!\!\sum_{\mathbf b\in [m^k]^s\setminus\{\mathbf a\}} \!\!\underline\lim_{\,n\to\infty} p_n(\mathbf b)\leqslant 1-\!\!\!\sum_{\mathbf b\in [m^k]^s\setminus\{\mathbf a\}}\!\!\! V(J_k(\mathbf b))=V(J_k(\mathbf a))
$$
(in the first passage to the limit we use the fact that the sum of lower limits does not exceed the lower limit of the sum).

Thus,\quad $\lim\limits_{\,n\to\infty} p_n(\mathbf a)=V(J_k(\mathbf a))$.
$\square$

In other words, we have proved that the following condition is sufficient for the uniformity with each suffix: for any collection of cubes $J_k(\mathbf a)$ the probability of hitting the union of their neighborhoods $O_k(\mathbf a)$ should be asymptotically (with $n$ tending to infinity) bounded from below by the total volume of considered cubes.

The following result is one of simple applications of propositions proved in this section.
\begin{lm}
\label{const}
If a collection of functions $(f_1,\ldots,f_s)$ is uniform with each suffix, then so is the collection $(f_1+z_1,\ldots,f_s+z_s)$, where $z_1,\ldots,z_s$ are arbitrary integer constants.
\end{lm}
\par\noindent\textbf{Proof:}
Let $\ell$ be the maximal number of digits in the $m$-ary notation that are necessary for writing $|z_1|,\ldots,|z_s|$; $k\in\mathbb N$, $n\geqslant \ell+k$, $x\in [m^n]$, and $y_i=substr(f_i(x), n, n-k)$, $i=1,\ldots,s$. By definition, ${substr(|z_i|,n,n-k)=0}$. Consequently, by Lemma~\ref{stolbik}, 
\begin{equation}
\label{zz}
(substr(f_i(x)\pm |z_i|, n, n-k) - y_i) \bmod m^k \in \{0, 1,m^k-1\}.
\end{equation}
In other words, the replacement of $f_i$ by $f_i+z_i$ makes elements of the vector $substr((f_1 (x),\ldots,f_s(x)),n,n-k)$ change at most by 1.

If a collection of functions $(f_1,\ldots,f_s)$ is uniform with each suffix, then condition~(\ref{mainmaindef}) is fulfilled. This condition means that the inequality 
$$
{\mathbb P}(\varphi_n(f_1,\ldots,f_s)\in J_k(\mathbf a)\,|\, \xi_d=\beta)\geqslant 
V(J_k(\mathbf a))-\varepsilon
$$
is valid for any arbitrarily small $\varepsilon$ and sufficiently large~$n$. Consequently, for any set $M\subseteq [m^k]^s$ with sufficiently large~$n$, 
$$
{\mathbb P}(\varphi_n(f_1,\ldots,f_s)\in \bigcup_{\mathbf a\in M} J_k(\mathbf a)\,|\, \xi_d=\beta)\geqslant  
\sum_{\mathbf a\in M} V(J_k(\mathbf a))-|M|\varepsilon.
$$

Note that if for an elementary event $\omega\in\Omega$, 
$$
\varphi_n(f_1,\ldots,f_s)\in \bigcup_{\mathbf a\in M} J_k(\mathbf a),
$$
then, as follows from~(\ref{zz}), for the same elementary event we get 
$$
\varphi_n(f_1+z_1,\ldots,f_s+z_s)
\in \bigcup_{\mathbf a\in M} \bigcup_{\mathbf c\in O_k(\mathbf a)}J_k(\mathbf c).
$$
Evidently, the same implication is valid even with an additional constraint $\xi_d=\beta$ imposed on the set of elementary events. Therefore, for sufficiently large~$n$,
$$
{\mathbb P}(\varphi_n(f_1+z_1,\ldots,f_s+z_s)
\in \bigcup_{\mathbf a\in M} \bigcup_{\mathbf c\in O_k(\mathbf a)}J_k(\mathbf c)\,|\, \xi_d=\beta)\geqslant 
\sum_{\mathbf a\in M} V(J_k(\mathbf a))-|M|\varepsilon.
$$
Thus, functions $f_1+z_1,\ldots,f_s+z_s$ satisfy all conditions of Lemma~\ref{OK}.
$\square$

\section{The Weyl criterion and linear combinations}
\label{weil_section}
Having proved Theorem~\ref{monoms}, let us prove that linearly independent polynomials with integer coefficients whose free term equals 0 are also uniform with each suffix. Evidently, this proposition follows from the next lemma; {since we also use its statement in the induction process, we give it here.}

\begin{lm}
\label{B}
Let $A$ be a nondegenerate $s\times s$-matrix with integer elements. If a collection of functions $f=(f_1,\ldots,f_s)$ generates a uniform distribution or it is uniform with each suffix, then so is the collection $g=(g_1,\ldots,g_s)$, where
\begin{equation}
\label{gf}
g^T=A f^T
\end{equation}
(here $T$ is the transposition sign; it means that in this case each row turns into a column).
\end{lm}

Recall that a discrete uniform distribution in a finite set $Y$ is defined by us in a standard way, namely, all elements of this set have equal probabilities ($1/|Y|$). In this section, we consider a discrete uniform distribution in a finite \textit{multiset}~$X$, whose elements can be repeated. In this case we understand a discrete uniform distribution as a discrete distribution on the set~$Y$ consisting of distinct elements~$y$ of the set~$X$ which is defined by probabilities
$$
{\mathbb P}(y)=\frac{|\{x: x=y, x\in X\}|}{|X|}.
$$
Thus, for example, for the discrete uniform distribution on the multiset $\{1,0,1\}$ it holds that ${\mathbb P}(1)=2/3$, ${\mathbb P}(0)=1/3$.

\begin{utv}[The Weyl criterion]
\label{two}
Let $\zeta_n$ ($n=1,2,\ldots$) be random vectors with the discrete uniform distribution in finite multisets $\mathbf X(n)$, $\mathbf X(n)\in\mathbb R^s$. Denote by $(\cdot\bmod 1)$ the operation of extracting a fractional part of each component of an $s$-dimensional vector and do by $\langle \cdot,\cdot \rangle$ the sum of componentwise products of $s$-dimensional vectors. The sequence $\zeta_n\bmod 1$ weakly converges to the continuous uniform distribution on $[0,1)^s$ if and only if for any $\mathbf h\in \mathbb Z^s$, $\mathbf h\ne 0$,
\begin{equation}
\label{equtvtwo}
\lim_{n\to\infty} \frac1{|\mathbf X(n)|}\sum\limits_{\mathbf x\in \mathbf X(n)} \exp(2\pi i \langle \mathbf{h},\mathbf{x} \rangle)=0.
\end{equation}
\end{utv}

For completeness of the study, let us prove this proposition here (in fact we somewhat modify the proof given in \cite[sections~1.1 and~1.2]{nidera} for our case). 
\par\noindent\textbf{Proof:}
Note that the expression under the limit sign in formula~(\ref{equtvtwo}) represents the mean of the random value $\exp(2\pi i \langle \mathbf{h}, \zeta_n\rangle)$, i.e.,
\begin{equation}
\label{char}
\frac1{|\mathbf X(n)|}\sum\limits_{\mathbf x\in \mathbf X(n)} \exp(2\pi i \langle \mathbf{h},\mathbf{x} \rangle)=
\mathbb E\exp(2\pi i \langle \mathbf{h}, \zeta_n\rangle).
\end{equation}
Evidently, equality~(\ref{char}) remains valid even with the random vector $\zeta_n\mod 1$ in place of $\zeta_n$ in the right-hand side of equality~(\ref{char}).

Denote by $\zeta$ the random vector with the continuous uniform distribution on $[0,1)^s$. One can easily make sure that for any $\mathbf h\in\mathbb Z^s$, $\mathbf h\neq 0$,
$$
\mathbb E\exp(2\pi i \langle \mathbf{h}, \zeta\rangle)=0.
$$
In addition, it is evident that with $\mathbf h=0$ for any $x\in\mathbb R^s$ it holds that $\exp(2\pi i \langle \mathbf{h}, \mathbf{x}\rangle)\equiv 1$. Therefore with $\mathbf h=0$ for any $s$-dimensional random vector $\kappa$,
$$
\mathbb E\exp(2\pi i \langle \mathbf{h}, \kappa\rangle)=1.
$$
Therefore, assumptions of Proposition~\ref{two} mean that for any $\mathbf{h}\in\mathbb Z^s$,
\begin{equation}
\label{Veillimit}
\lim_{n\to\infty}\mathbb E\exp(2\pi i \langle \mathbf{h}, \zeta_n\bmod 1\rangle)=
\mathbb E\exp(2\pi i \langle \mathbf{h}, \zeta\rangle).
\end{equation}
We have to prove that these equalities are necessary and sufficient conditions for a weak convergence of $\zeta_n\bmod 1$ to $\zeta$. 

As is well known (see, for example, \cite{billing}), the sequence of $s$-dimensional random vectors $\kappa_n$ weakly converges to the random vector $\kappa$ if and only if for any continuous bounded complex-valued function~$\Psi$ ($\Psi$: $\mathbb{R}^s\to \mathbb{C}$) it holds that
\begin{equation}
\label{wiki}
\lim_{n\to\infty}\mathbb E \Psi(\kappa_n)=\mathbb E \Psi(\kappa).
\end{equation}
Note that without loss of generality we can consider only real-valued functions~$\Psi$, this form of equality~(\ref{wiki}) is often used as a definition of the weak convergence notion. 

Evidently, if $\kappa_n=\zeta_n\bmod 1$ and $\kappa=\zeta$, then equality~(\ref{Veillimit}) represents a particular case of equality~(\ref{wiki}), therefore correlations~(\ref{Veillimit}) with $\mathbf h\in\mathbb{Z}^s$ are necessary conditions for the weak convergence.

Let us prove the sufficiency. Since $\kappa_n$ take on values in $[0,1)^s$, it suffices to verify correlation~(\ref{wiki}) for continuous bounded complex-valued functions $\Psi$ defined on $[0,1)^s$. It is convenient to extend this function up to a periodic one on $\mathbb R^s$ by putting $\Psi(\mathbf x)=\Psi(\mathbf x\bmod 1)$. If this extension is also continuous, then we treat the function~$\Psi$ as \textit{continuous on a torus}. Let us first prove equality~(\ref{wiki}) for this case.

Let $\varepsilon$ be an arbitrary positive number. In the case under consideration, according to the Weierstrass theorem, there exists a trigonometric polynomial $\Psi_0(\mathbf x)$, i.e., a finite linear combination of functions $\exp(2\pi i \langle \mathbf{h}, \mathbf x\rangle)$, $\mathbf{h}\in\mathbb Z^s$, with complex coefficients such that
\begin{equation}
\label{sup}
\sup_{\mathbf x \in [0,1)^s}|\Psi_0(\mathbf x)-\Psi(\mathbf x)|\leqslant \varepsilon.
\end{equation}
We get
$$
|\mathbb E (\Psi(\zeta)-\Psi(\zeta_n))|\leqslant |\mathbb E (\Psi(\zeta)-\Psi_0(\zeta))|+
|\mathbb E (\Psi_0(\zeta)-\Psi_0(\zeta_n))|+|\mathbb E (\Psi_0(\zeta_n)-\Psi(\zeta_n)|.
$$
In view of~(\ref{sup}) the first and third terms in the right-hand side of this correlation do not exceed $\varepsilon$ for any~$n$. In view of~(\ref{Veillimit}) with sufficiently large~$n$ the second term does not exceed $\varepsilon$. Therefore, for functions $\Psi$, which are continuous on a torus, equality~(\ref{wiki}) is proved.

Furthermore, it is evident that for any continuous and bounded on the cube $[0,1)^s$ function $\Psi(\mathbf x)$ there exist two functions $g_1(\mathbf x)$ and $g_2(\mathbf x)$ which are continuous on a torus (extendable to functions which are continuous and periodic on $\mathbb R^s$) and such that for all $\mathbf x\in [0,1)^s$,
\begin{equation}
\label{g1g2}
g_1(\mathbf x)\leqslant\Psi(\mathbf x)\leqslant g_2(\mathbf x)
\end{equation}
and $\int_{\mathbf x\in [0,1)^s} (g_2(\mathbf x)-g_1(\mathbf x))\,d\mathbf x \leqslant \varepsilon$
(cf. with~\cite[p.~3]{nidera}). The latter inequality means that
$$
\mathbb E g_2(\zeta)-\varepsilon\leqslant \mathbb E g_1(\zeta),
$$
inequalities~(\ref{g1g2}) can also be written in terms of mean values with random vectors in place of $\mathbf x$. We get
$$
\mathbb E\Psi(\zeta)-\varepsilon\leqslant\mathbb E g_2(\zeta)-\varepsilon\leqslant\mathbb E g_1(\zeta)=
\lim_{n\to\infty}\mathbb E g_1(\zeta_n)\leqslant \underline\lim_{\,n\to\infty}\mathbb E\Psi(\zeta_n).
$$
One can analogously prove that $\overline\lim_{\,n\to\infty}\mathbb E\Psi(\zeta_n)\leqslant \mathbb E\Psi(\zeta)+\varepsilon$, namely,
$$
\overline\lim_{\,n\to\infty}\mathbb E\Psi(\zeta_n)\leqslant
\lim_{n\to\infty}\mathbb E g_2(\zeta_n)=\mathbb E g_2(\zeta)\leqslant
\mathbb E g_1(\zeta)+\varepsilon\leqslant\mathbb E\Psi(\zeta)+\varepsilon.
$$
Due to the arbitrariness of the choice of $\varepsilon$ we conclude that $\lim_{n\to\infty}\mathbb E \Psi(\zeta_n)=\mathbb E \Psi(\zeta)$.
$\square$

One can easily prove the next proposition with the help of the Weyl criterion.

\begin{sld}
\label{sl:two}
Let assumptions of Proposition~\ref{two} be fulfilled and let a sequence $\zeta_n\bmod 1$ weakly converge to the continuous uniform distribution on $[0,1)^s$. Then the sequence $\zeta'_n\bmod 1$, where
$\zeta'_n=A \zeta_n$, while $A$ is a nondegenerate $s\times s$-matrix with integer elements, also weakly converges to the continuous uniform distribution on $[0,1)^s$.
\end{sld}
\par\noindent\textbf{Proof:}
Really, with $\mathbf h\in \mathbb Z^s$, $\mathbf h\ne 0$, we get
$$
\frac1{|\mathbf X(n)|}\sum\limits_{\mathbf x\in \mathbf X(n)} \exp(2\pi i \langle \mathbf{h},A \mathbf{x} \rangle)=
\frac1{|\mathbf X(n)|}\sum\limits_{\mathbf x\in \mathbf X(n)} \exp(2\pi i \langle A^T \mathbf{h},\mathbf{x} \rangle).
$$
Since components of the vector $A^T \mathbf{h}$ are integer and (due to the nondegeneracy of the matrix~$A$) nonzero, the latter sequence tends to $0$ as $n$ tends to infinity.
$\square$

\par\noindent\textbf{Proof of Lemma~\ref{B}.}
Let us define the operation of extracting the fractional part $(\cdot\bmod 1)$ for an arbitrary $m$-adic vector (whose elements belong to the set~$\mathbb Q_m$); in~\cite{Gelfand} this operation is denoted by the symbol $\{\cdot\}$. Namely, for any~$m$-adic number~$x$ in form~(\ref{series}) we put
$$
x\bmod 1=\sum_{i=k}^{-1} a_i m^i.
$$
In the case of a vector $\mathbf x$, $\mathbf x\in \mathbb Q_m^s$, the operation $\mathbf x\bmod 1$ is performed componentwisely.

Note that for rational numbers in form $\ell/m^k$, where $\ell\in\mathbb Z$, $k\in\mathbb N$, the operation $(\cdot\bmod 1)$ (understood both in the usual sense and in the $m$-adic one) gives one and the same result. Moreover, for any $a\in\mathbb Z$ and $x\in \mathbb Q_m$ we conclude that
$$
(a x)\bmod 1=(a (x\bmod 1))\bmod 1,
$$
because the multiplication by an integer number does not exceed the length of the fractional part. Analogously, for $\mathbf x\in \mathbb Q_m^s$ and for any matrix~$A$ with integer elements we get
\begin{equation}
\label{Amathbbx}
(A \mathbf x)\bmod 1=(A (\mathbf x\bmod 1))\bmod 1.
\end{equation}

In view of the 1-Lipschitz property of considered functions~$f_i$ for any elementary event $\omega\in \mathbb Z_m$, formula~(\ref{otherreptresentation}) for calculating the vector $\varphi_n (f_1,\ldots,f_s)$ gives $\varphi_n (f_1,\ldots,f_s)=((f_1(\omega),\ldots,f_s(\omega))/{m^n})\bmod 1$. Let collections of functions $g$ and $f$ obey correlation~(\ref{gf}). Applying~(\ref{Amathbbx}), we conclude that
\begin{gather}
\varphi_n (g_1,\ldots,g_s)^T=((g_1(\omega),\ldots,g_s(\omega))^T/{m^n}) \bmod 1=\notag\\
=(A\ (f_1(\omega),\ldots,f_s(\omega))^T/{m^n}) \bmod 1=(A\ \varphi_n (f_1,\ldots,f_s)^T)\bmod 1.\notag
\end{gather}
This allows us to immediately apply Corollary~\ref{sl:two}.

In the case when a collection of functions $f=(f_1,\ldots,f_s)$ immediately generates the uniform distribution, we use $\varphi_n (f_1(\omega),\ldots,f_s(\omega))^T$ for $\zeta_n$. In this case the set $\textbf X(n)$ has the cardinality number $m^n$ and represents the multiset
\begin{equation}
\label{Xn}
\left\{\left(\frac{(f_1 \bmod m^n)(x)}{m^n},\ldots,\frac{(f_s\bmod m^n)(x)}{m^n}\right): x\in [m^n] \right\}.
\end{equation}
Therefore, the assertion of the lemma that if a collection of functions $f=(f_1,\ldots,f_s)$ generates a uniform distribution, then so does the collection $g=(g_1,\ldots,g_s)$ follows from Corollary~\ref{sl:two}.

The case when a collection of function $f=(f_1,\ldots,f_s)$ is uniform with each suffix~$\beta=(\beta_1,\ldots,\beta_d)$ can be studied analogously. Therefore the cardinality number of the set $\textbf X(n)$ equals $m^{n-d}$; the set $\textbf X(n)$ represents only a part of the multiset~(\ref{Xn})(rather than the entire one), namely, it depends only on $x\in [m^n]$ which have a fixed suffix~$\beta$ (rather than on all of them).
$\square$

\section{The proof of the uniformity of the collection of monomials}
\label{monoms_section}

In this section we prove Theorem~\ref{monoms}.
\par\noindent\textbf{Proof:} Let us prove the theorem by induction with respect to~$s$.

The induction base for $s=1$ is evident, namely, the condition ensuring that $\frac{x}{m^n}$ belongs to the semiinterval $J_k(a_1)$ is equivalent to that $substr(x,n,n-k)=a$, therefore, the desired correlation~(\ref{mainmaindef}) is fulfilled as the equality for any $n \geqslant d+k$.

According to the induction principle, it suffices to prove the fulfillment of assumptions of Lemma~\ref{OK} for the collection $x,x^2,\ldots,x^{s}$ for each fixed suffix~$\beta$, $\beta\in[m^d]$. This means that for any arbitrarily small $\varepsilon$ it suffices to find (for each cube $J_K(\mathbf b)$, $\mathbf b\in [m^K]^s$) a collection of elementary events $\mathcal A_{d,K}(\beta,\mathbf b,n)$ such that the following conditions are fulfilled with sufficiently large~$n$:
\begin{enumerate}
\item[A.] For any $\omega\in \mathcal A_{d,K}(\beta,\mathbf b,n)$ it holds that
\begin{enumerate}
\item[$\alpha$)] $\omega\bmod m^d=\beta$
\item[$\beta$)]  $\varphi_n(x,x^2,\ldots,x^{s})(\omega)\in \bigcup_{\mathbf c\in O_K(\mathbf b)} J_K(\mathbf c)$;
\end{enumerate}
\item[B.] ${\mathbb P} (\mathcal A_{d,K}(\beta,\mathbf b,n))\,m^{-d}  \geqslant V(J_K(\mathbf b))-\varepsilon$;
\item[C.] $\mathcal A_{d,K}(\beta,\mathbf b,n)\cap \mathcal A_{d,K}(\beta,\mathbf b',n)=\emptyset$, if $\mathbf b\neq \mathbf b'$, $\{\mathbf b,\mathbf b'\} \subseteq [m^K]^s$.
\end{enumerate}
The latter requirement ensures that assumptions of Lemma~\ref{OK} are fulfilled not only for $M=\{\mathbf b\}$, but also for the case when the set~$M$ consists of several elements.

Note that elementary events $\omega$ belong to $\mathbb Z_m$, and conditions imposed by us are connected only with random values $\xi_n$; $\xi_n=\omega\bmod m^n$ (which are discretely uniformly distributed on the set $[m^n]$). We use the trivial correspondence between sets $\mathcal A_{d,K}(\beta,\mathbf b,n)$ and $\mathcal A'_{d,K}(\beta,\mathbf b,n)=\{x:
x=\omega\bmod m^n, \omega\in \mathcal A_{d,K}(\beta,\mathbf b,n)\}$, the first of them are obtained from the latter ones by adding all possible values from $m^n \mathbb Z_m$. It is more convenient to operate with terms of sets $\mathcal A'_{d,K}(\beta,\mathbf b,n)$, all whose elements belong to~$[m^n]$.
Requirement~B means that {the ratio of the measure} of
$\mathcal A'_{d,K}(\beta,\mathbf b,n)$ {to the measure of} all $x\in [m^n]$, $x\bmod m^d =\beta$, is bounded from below by the value $V(J_K(\mathbf b))-\varepsilon$, i.e., $m^{-sK}-\varepsilon$.

The induction hypothesis is that the collection $x,x^2,\ldots,x^{s-1}$ is uniform with each suffix. Then by Lemma~\ref{B} so is the collection $2 x,\,3 x^2,\ldots,s x^{s-1}$. Therefore, we can apply Corollary~\ref{sldcond} to it and assume that for $(f_1,\ldots,f_{s-1})=(2 x,\,3 x^2,\ldots,s x^{s-1})$ there exists~$N$ such that condition~(\ref{cond}) is fulfilled.

Fix parameters $j$ and $k$ in Corollary~\ref{sldcond}, namely,
$$
j=s K,\qquad k=2(s-1)K,
$$
and consider the collection $a_i=m^{(s-1)K+(i-1)K}$, $i=1,\ldots,{s-1}$, for $\mathbf a$. Therefore, in the $m$-ary notation numbers $a_i$ have only one unit (each one is located on its own place), while the rest digits are zeros. Let us prove that then with $n\geqslant 2N$ (where $N$ is such that condition~(\ref{cond}) is fulfilled) we can find sets $\mathcal A'_{d,K}(\beta,\mathbf b,n)$ with the desired properties.

In what follows, for convenience, we assume that the value~$\ell$ in condition~(\ref{cond}) is calculated in some deterministic way, for example, as the least value among those that satisfy this condition.

Let us first define~$\mathcal B \subseteq [m^n]$ as the totality of all~$z\in [m^n]$ such that 1) $x=substr(z,N,0)$ satisfies requirements of condition~(\ref{cond}), 2) for calculated~$\ell$ it holds that $substr(z,n-\ell+k/2,n-\ell)=0$, and, finally, 3) $substr(z,n,n-K)=0$. Note that the latter two constraints are independent, because $\ell\geqslant j=sK$. These constraints ensure that exactly $sK$ digits of $z$ written in the $m$-ary notation are equal to $0$. Evidently, the probability of this event (when all values of these digits have equal probabilities) equals $m^{-sK}$.

The {first} constraint is also independent of the {last}
two ones ($\ell\leqslant N \leqslant n/2$). Since the relative number of~$x$ satisfying the first constraint (among all~$x$ such that $x\bmod m^d=\beta$) exceeds $1-\varepsilon$, the relative number of $z$ satisfying all three constraints (among all~$z$ such that $z\bmod m^d=\beta$) is greater than
\begin{equation}
\label{eneq}
(1-\varepsilon) m^{-s K}> m^{-s K}-\varepsilon.
\end{equation}

Let us construct an invertible operator~$L_z$, $z\in \mathcal B$, that maps collections $\mathbf b\in [m^{K}]^s$ to those $\mathbf c=(c_0, c_1,\ldots,c_{s-1})\in [m^{K}]^s$ and has the following property. Put
\begin{equation}
\label{lll}
y=z+c_0 m^{n-K}+c\, m^{n-\ell},\qquad \mbox{where $c=\sum_{i=1}^{s-1} c_{s-i}\, m^{K(i-1)}$,}
\end{equation}
here the parameter $\ell$ is defined by condition~(\ref{cond}) for $x=z\bmod p^N$. Then
\begin{equation}
\label{last}
\varphi_n(y, y^2,\ldots, y^s) \in \bigcup_{\mathbf b' \in O_K(\mathbf b)} J_K(\mathbf b').
\end{equation}

Note that $c\in [m^{k/2}]$, because $c_i\in [m^K]$. Therefore, $y$ (in the $m$-ary notation)
is obtained from $z$ by placing $c_i$ instead of $sK$ zeros
{whose positions are defined by} conditions~2) and~3). Therefore (in view of~(\ref{eneq})) the totality of all such~$y$ for fixed~$\mathbf b$ forms a set $\mathcal A'_{d,K}(\beta,\mathbf b,n)$ which satisfies requirements A--C.

Let us construct the operator~$L_z$, as was stated earlier.
Put $c_0=b_1$; $z_1=z+c_0 m^{n-K}$, where $z\in \mathcal B$, and $x$,
which is defined above, coincides with $z_1\bmod m^N=y\bmod m^N$.
Note that the value $y$ is representable as $z_1+c\, m^{n-\ell}$.

Since $c_0 = substr(z_1, n, n-K)=substr(y, n, n-K)$,
we conclude that $y/m^n\in J_{K}(b_1)$.
Furthermore, let us consider $y^t\bmod m^n$, $t=2,\ldots,s$.

Let us apply the binomial theorem to this expression. Since we perform all calculations modulo $m^n$, it suffices to restrict ourselves to the first two terms. Thus,
$$
y^t\bmod m^n=(z_1^t\bmod m^n+t z_1^{t-1}\,c\, m^{n-\ell}\bmod m^n)\bmod m^n.
$$
The latter term in the inner sum is representable as
$$
((t z_1^{t-1}\bmod{m^\ell})\, c\, m^{n-\ell})\bmod m^n=
((t x^{t-1}\bmod m^\ell)\,c\, m^{n-\ell})\bmod m^n.
$$
When deducing the latter equality, we have taken into account the fact that the operation of raising to a power is an {1-Lipschitz} function.

According to Corollary~\ref{sldcond}, the number $t x^{t-1}\bmod m^\ell$ is representable as $x'+a_t m^{\ell-k}$, where $x'\in [m^{\ell-k}]$ and $a_t=m^{(s-1)K+(t-2)K}$, while $t=2,\ldots,s$. Let us now multiply this sum by $c\, m^{n-\ell}$. To this end, let us separately calculate the product of each term. 

The multiplication of the second term by $c\, m^{n-\ell}$ gives $a_t m^{\ell-k} c\, m^{n-\ell}=c\, m^{n-K-(s-t)K}$. Furthermore, $x'\in [m^{\ell-k}]$, $c\in[m^{k/2}]$, consequently, $x' c\in [m^{\ell-k/2}]$,\ $x' c\, m^{n-\ell}\in [m^{n-k/2}]$. As a result,
$$
substr((t z_1^{t-1})\,c\, m^{n-\ell},n,n-K)=substr(c\, m^{n-K-(s-t)K},n,n-K)=c_{t-1}.
$$
By Lemma~\ref{stolbik} we conclude that $substr(y^t,n,n-K)$ differs from $$(c_{t-1}+substr(z_1^t,n,n-K))\bmod m^K$$ at most by 1.

This allows us to (rather easily) define the operator~$L_z:\mathbf c=L_z(\mathbf b )$ which satisfies~(\ref{last}). Put $c_0=b_1$, $z_1=z+c_0 m^{n-K}$, $$c_{t-1}=(b_t-substr(z_1^t,n,n-K))\bmod m^K,\quad  t=2,\ldots,s.$$ Evidently, this map represents a bijection.

Define the set $\mathcal A'_{d,K}(\beta,\mathbf b,n)$ as the totality of all $y$ given by formula~(\ref{lll}) for all $z\in \mathcal B$. By construction, the set $\mathcal A'_{d,K}(\beta, n, \mathbf b)$ is the desired one, whence it follows that the set $\mathcal A_{d,K}(\beta, n, \mathbf b)$ satisfies requirements A--C.
$\square$

\section{The uniformity in the case of polynomials and the proof of the main theorem}
\label{polynoms_section}
In this section, using previous results, we easily prove Theorem~\ref{mainthm}. It represents a corollary of the following proposition.

\begin{thm}
\label{th:last}
For any polynomial $f$ with integer coefficients, whose power exceeds 1, and any $s\in\mathbb N$, the collection of functions $x,f,f^{(2)},\ldots,f^{(s-1)}$ (where $x$ is the identical map) is uniform with each suffix.
\end{thm}

Let us first state several simple but important propositions; one can easily deduce them from previously proved ones.
\begin{thm}
\label{polys}
Let $A$ be an arbitrary nondegenerate $s\times s$-matrix with integer elements. Define a column of polynomials $f_1,\ldots,f_s$ by the correlation $(f_1,\ldots,f_s)^T=A (x,\ldots,x^{s})^T + z$, where $z$ is an arbitrary constant $1\times s$-column, all whose elements are integer. Then the collection of functions $f_1,f_2, \ldots,f_s$ is uniform with each suffix.
\end{thm}
\par\noindent\textbf{Proof:}
The assertion of this theorem evidently follows from lemmas~\ref{B} and~\ref{const} (in view of Theorem~\ref{monoms}).
$\square$

\begin{lm}
Let a collection of functions $f_1,f_2, \ldots,f_s$ be uniform with each suffix. Then any subcollection $f_{i_1},f_{i_2},\ldots,f_{i_k}, 1 \leqslant i_1 < i_2 < \ldots < i_k \leqslant s$ is also uniform with each suffix.
\end{lm}
\par\noindent\textbf{Proof:} 
By the definition of the conditional probability, if events $A_1,\ldots A_N$ do not intersect, then
\begin{equation}
\label{eqlast}
\sum_{i=1}^N \mathbb P(A_i|B)=\mathbb P(\cup_{i=1}^N A_i|B).
\end{equation}

Let us perform the summation in the left- and right-hand sides of formula~(\ref{mainmaindef}) with respect to all possible $a_i\in [m^k], 1 \leqslant i \leqslant n, i \notin \{i_1,i_2,\ldots,i_k\}$ and interchange the sum and limit signs in the left-hand side. Using formula~(\ref{eqlast}), we get~(\ref{mainmaindef}) for the subcollection $f_{i_1},f_{i_2},\ldots,f_{i_k}$.
$\square$

\par\noindent\textbf{Proof of Theorem~\ref{th:last}:} Let the power of $f^{(s-1)}$ equal $d$. Evidently, the collection of polynomials $x,f, f^{(2)}, \ldots,f^{(s-1)}$ contains no polynomials of one and the same degree. Let us complement it with arbitrary polynomials such that for each $i,1\leqslant i \leqslant d$ the resulting set contains exactly one polynomial of degree $i$. According to Theorem~\ref{polys} (since any triangular matrix is nondegenerate), this collection of functions is uniform with each suffix. Since $x,f, f^{(2)}, \ldots,f^{(s-1)}$ is its subcollection, in view of the previous lemma this collection is also uniform with each suffix, which was to be proved.
$\square$

\par\noindent\textbf{Proof of theorems~\ref{mainthm} (and~\ref{mainthm0}):}
In accordance with Proposition~\ref{utv_ochevidno} and Theorem~\ref{th:last}, the collection of functions $x, f, f^{(2)}, \ldots, f^{(s-1)}$ generates a uniform distribution. Therefore, Theorem~\ref{mainthm} is valid. 
In view of Proposition~\ref{ochevidno2} this proves Theorem~\ref{mainthm0}.
$\square$

\section{Conclusion}

In this paper we prove that the sequence constructed by iterations of an arbitrary polynomial with integer coefficients, whose degree is at least 2 modulo $m^n$, generates a uniform distribution as $n$ tends to infinity.

Moreover, we also prove a more general assertion, namely, we prove that any collection of polynomials with integer coefficients generates a uniform distribution modulo $m^n$ as $n$ tends to infinity for any~$m$, provided that these polynomials (without free terms) are linearly independent.

The discrepancy bound {that can be established with the help of the technique used in this paper} is weaker than that established in~\cite{nemcy},\cite{nemcy2},\cite{nem3} for concrete classes of polynomials. Recall that we understand a discrepancy $D_{m^n}(\mathbf P_n^s(f))$ as 
$$\sup_{J \in \mathfrak J_s} \left|\frac{|\mathbf P_n^s(f) \cap J|}{|\mathbf P_n^s(f)|} - V(J)\right|$$
(here we use denotations given at the beginning of Section~2). The bounds established in papers~\cite{nemcy},\cite{nemcy2},\cite{nem3} give the main term of the asymptotics $D_{m^n}$ in the form $m^{nc}$, where $c=-1/2$, with some logarithmic corrections consistent with the law of the iterated logarithm (note that there exist polynomials which violate this bound, see~\cite{nem3}). The discrepancy bound that can be established with the help of the technique used in this paper allows us only to prove that the upper bound for $D_{m^n}$ decreases as $\log n$ raised to some negative power.

Moreover, it is possible to establish a uniformity criterion for a collection of functions for a finite automaton. As appeared, an automaton maps a collection of functions generating a uniform distribution to some other collection with the same property if and only if the synchronization takes place (see, for example,~\cite{volkov}). Therefore, instead of linear combinations (see~Lemma~\ref{B}) one can use a more general structure~\cite{dan,CIAA}.

The author is grateful to his scientific supervisor V.S.~Anashin for the problem statement and for useful discussions.
The author also is grateful to M.A.~Cherepnev for the comments.

\end{document}